\documentclass[11pt]{article}
\usepackage{amsmath,amssymb,amsfonts,latexsym}
\usepackage{amssymb}
\usepackage{amsthm}
\usepackage{cancel}
\usepackage{newlfont}
\usepackage{youngtab}

\makeatletter
\def \fixequation{\let \c@equation\c@theorem\let \p@equation\p@theorem\let \theequation \thetheorem}
\makeatother \theoremstyle{plain}
\newtheorem{theorem}{Theorem}
\newtheorem{lemma}{Lemma}

\theoremstyle{definition}

\theoremstyle{remark}

\begin{document}
\def\s{\sigma}
\def\A{\mathcal{A}}
\def\B{\mathcal{B}}
\def\ga{\gamma}
\def\La{\Lambda}
\def\la{\lambda}
\def\p{\phi}
\def\ep{\epsilon}
\def\sg{\sigma}
\def\al{\alpha}
\def\gp{\frak P}
\def\gq{\frak q}
\def\gz{\frak Z}
\def\gs{\frak S}
\def\gg{\frak g}
\def\kin{k \in [1,n]}
\def\kim{k \in [1,m]}
\def\ss{s^k_1, s^k_2}
\def\Jin{J \subset [1,n]}
\def\ts{\tilde{s}}
\def\mN{\mathbb{N}^+}
\def\mz{\mathbb{Z}}
\def\mp{\mathbb{P}}
\def\mzt{\mathbb{Z}_2}
\def\Cb{\mathbf{C}}
\def\bE{\mathbf{E}}
\def\qs{q_{\mbox{SAT}}}
\def\sbe{\mbox{SAT}(\mathbf{E})}
\def\rx{R[X]}
\def\rxo{R[X_0]}
\def\on{\otimes^N \mathbb{C}^n}
\def\onn{\otimes^n \mathbb{C}^n}
\def\one{\otimes^N \mathbb{C}^n(\eta)}
\def\noo{(v_{i_1} \otimes \cdots \otimes v_{i_N})}
\def\GL{GL_n(\mathbb{C}}
\def\sn{S_N(U_n)}
\def\mbc{\mathbb{C}}
\begin{center} {\bf \Large Averaging Over the Unitarian Group and the Monotonicity Conjecture of Merris and Watkins \\by Avital Frumkin }
\end{center}
\large {abstruct}\\
 We show that the monotonicity conjecture of
Merriss and Watkins is true in average when taking the set of
matrices of given non negative  spectra as probability space with
respect to the Haar measure of the unitarian group .

\section{Introduction}
Let $\mbc S_n$ be the group algebra of the symmetric group $S_n$
over the complex numbers.   ,$\mbc S_n$ includes all the functions
from $S_n$ to $\mbc$ as vector space of dimension N! . The
multiplication is defined as convolution of group's functions , i.e.
\begin{equation}
f \cdot g(\sigma) = \sum_{\gamma \in S_n} f(\gamma) \ g(\gamma^{-1} \sigma)
\end{equation}
Given $n \times n$ complex matrix $A$, we define a function from
$\mbc S_n$ to $\mbc$ by the formula
\begin{equation}
f \in \mbc S_n \quad \quad f \rightarrow d_f (A) = \sum_{\sigma \in
S_n} f(\sigma)\prod_i A_{i \ \sigma(i)}
\end{equation}   ( d is for "determinant")\\
In these notations, the determinant of a matrix $A$ is given by the
formula
\begin{equation} det A = d_f (A) ,  f : S_n \rightarrow \mbc , f(\sigma)
= sig(\sigma).
\end{equation}
The permanent of $A$ is
\begin{equation}
d_f (A) ,   f : S_n \rightarrow \mbc, f(\sigma) = 1 . for  any
\sigma \in S_n.
\end{equation}

 For  given irreducible  character$\chi$ of $S_n$,after  identifying  it as a  function from $\mbc S_N$to the complexes   $d_\chi (A)$
is called the immanent of $A$ corresponding with $\chi$.

a partition of $N$ is a  vector of integers $\eta = \eta_1, \eta_2
\cdots$ so that $\eta_1 \geq \eta_2 \cdots \geq 0$ and $\sum \eta_1
=N$. We denote $\eta$ A  partition of $N$ by $\eta \vdash N$.

 The irreducible characters of $S_N$  correspond to
partitions of $N$, hence each irreducible character is attached to a
partition, .  So we denote $d_ \eta (A)$ instead of $d_{\chi_\eta}
(A)$.

The first theorem in the area of immanents is that of Schur in 1918
. One corollary of it is that

\begin{equation}
det (A) \leq \frac{d_\eta(A)}{\chi(e)}
\end{equation}
 to any $n\times n$ Hermitian positive
semi-definite matrix $A$ and $\eta \vdash n$
once $e$ is the trivial
element of $S_n$.

In my view,even the fact that $d_\eta(A)$ is non-negative to any
non-negative definite matrix $A$ and $\eta \vdash n$ is amazing.

 The word immanent, apparently invented by Littlewood,to
 these  complex functions  involving  $S_n$ group algebra and
$n \times n$ matrices over $\mbc$. See Chapter (6) in his book  The
Theory of Group Character [Li]. Littlewood includes in his
definition of immanent almost any monomials in the $n \times n$
matrices elements in his celebrated postulate about the
relationships between  immanents and Schur functions. Here is the
postulate ".... corresponding to any relation between Schur
functions of total weight $N$, we may replace the Schur function by
the corresponding immanents of complementary coaxial minor of
$[a_{st}]$ provided that every product is summed for all sets of
complementary coaxial minors."

This postulate, together with  Lieb's inequality [Li], is the corner
stone  of the paper of Merris and Watkins [M W].  In their paper,
they used it in innovative way  accomplishing inequalities and
identities of immanents and matrix functions  \\
Among other things, they raised a conjecture which was called
afterwards by Pate [P3] the Merris and Watkins (monotonicity)
conjecture . For treating  it, we need some definitions.\\

For a given partition of $N, \eta = \eta_1 \geq \eta_2 \geq \eta_n
\geq 0$. The  $\eta's$ weight space of $\on$ is spent  by all the
tensors $v_{i_1} \otimes \cdots \otimes v_{i_N}$ for which the
number of occurrences of $v_j$, i.e. the $j$ basis vector of
$\mbc^n$, is $\eta_j$.  The weight spaces in $\on$ are $S_N$ module
with respect to the Schur action of $S_N$ on $\on $:
 for $\sigma \in S_n$ its Schur action is defined by the formula\\
\begin{equation}
 \sigma (v_{i_1} \otimes \cdots \otimes v_{i_N} ) =
v_{i_\sigma(1)} \otimes \cdots \otimes v_{i_\sigma(N)}
 \end{equation}

. The $S_N$'s character corresponding with this  $\eta$ weight space
is denoted  $[\eta]$ in [M W]

In the language of character theory,  $[\eta] $ is the induce
character  from the trivial character of the Young sub group
corresponding with $\eta$ see [J K ]to $S_N$\\The dominant order on
the set of partitions of $N$  $\succ$,is defined by the formula
\begin{equation}
\eta \succ \eta' \Leftrightarrow \sum^j_{i=1} \eta_i \geq
\sum^j_{i=1} \eta'_i \quad \mbox{to any} \ j
\end{equation}

This order between partitions of $N$ has an intimate connection to the action of the Lie algebra $g \ell_n(\mbc)$ on $\on$
and through it with the action of  $GL_n(\mbc)$ on $\on$.

 By  using a clever  result in representation  theory of symmetric
groups ( theorem 22 of [JK])(among other results), Merris and
Watkins proved that for any non-negative definite $n \times n$
matrix $A$
\begin{equation}
\eta \succ \eta ' \Rightarrow d_{[\eta]} A \leq d_{[\eta ']} A
\end{equation}
On the other hand , they conjectured (God know how )that \\
\\
$\left( \begin{array}{l} n \\ \eta \end{array} \right)^{-1}
d_{[\eta]} A \geq \left( \begin{array}{l} n \\ \eta' \end{array}
\right)^{-1} d_{[\eta']} A$ . $\left( \begin{array}{c} n \\
\alpha  \end{array}\right)$ is the multinomial coefficient attached
to the partition $\alpha$.\\ Clearly $\left( \begin{array}{c} n \\
\alpha
\end{array}\right) $  is the dimension of the $\alpha 's$ weight
space in $\bigotimes^n\mbc$.

We shall prove (theorem 7)that when $\eta \succeq \eta'$ \\
\\
$\left(
\begin{array}{l} n \\ \eta \end{array} \right)^{-1} \int d_{[\eta']}
A^u du \geq \left( \begin{array}{l} n \\ \eta' \end{array}
\right)^{-1} \int d_{[\eta ']} A^u du$. when $A^u=uAu^{-1}$ and the
integration is done over the Unitarian group $U_n$ with respect to
the Haar measure  du of $U_n$ .

 When $n$ tends to infinity, and $\eta$ is a partition
of $n$ of a bounded number of parts, the character$[\eta]$

 tends to be closer to $\chi_\eta$, the irreducible
character corresponding with $\eta$, in a reasonable matric.More
precisely when $[\eta]$ is written as sum of ireduceabl characters
the contribution of those characters   far from $\chi_\eta$ to the
dimension of $[\eta]$ vanishes  in comparison to the full dimension.
See [F.G] or [M.H.Ch]( under the tittle  Kyel Werner theorem)

Because of  that and  the fact that the  induce characters are
easier than the irreducible one to compute, the asymptotic behavior
of them is worth attention when  looking for counter examples.Indeed
this is how I got to consider them . Hence before  turning  to the
Merris and Wotkins monotonicity  conjecture in average , we survey
the
monotonicity of immanents' status as I know it\\

 I believe  that the next theorem of
Pate  cover the  scope till now  \\
 Given a partition $\eta = \eta_1 \geq
\eta_2 \geq \cdots \eta_i
> \eta_{i+1} \cdots  $ let $ \eta' = \eta_1 \geq \eta_2 \geq \cdots
\eta_i -1\geq\eta_{i+1}\cdots\geq 1 \Rightarrow\frac{d_{\eta '}
(A)}{\chi_{\eta '} (e)} \leq \frac{d_{\eta} (A)}{\chi_\eta (e)}$.\\
 That  is to say that in removing   a corner in the Young diagram of shape
$\eta$ to the end decreases the normalized immanent corresponding
with the new partition.See [P2] , The first large  step in this
direction
was done in [P 1]\\

Our approach to computations of the $S_n$ matrix function begins
with an observation that appeared  in [Kos] where Kostant  reproved
that

$d_\eta (A) > d_\eta (I)$ for non-negative definite, or totally
positive matrices,A  of determinant 1 .\\
 His observation is that
$d_\eta(A)$ is the "trace" of $A^{\otimes n}$ when it acts on
$M_\eta(0)$.  The zero weight subspace of $M_\eta$ where $M_\eta$ is
the irreducible $GL_n(\mbc)$ module corresponding with $\eta$. The
zero weight space is the subspace of all the Tori invariant vectors.
In other words  the subspace   of the equipartition weight . The
parenthesis over the trace was needed since   $M_\eta(0)$ is not
respected by $A^{\otimes n}$ for a general matrix $A$ and
 some projection is needed. But by  averaging  the immanent of
$uAu^*$ over the Unitarian group with respect to its Haar measure ,
the projection, can be ignored  \\  Given $f\in \mbc S_N$ we denote
$\hat d_f(A)=\int d_f(A^u)du$ integrated over the unitarian
group with respect to the Haar measure \\
Let us denote   $\hat d_\eta (A)$ for $\hat d_{\chi_\eta} (A)$ . We
shall prove  (theorem 5)\\
\begin{equation}
\hat{d}_\eta(A) = trace \int A^{u \otimes n} du \left|_{M_{\eta
(0)}}\right . = \frac{s_\eta (A)}{s_\eta (I)} dim (\ M_\eta (0))
\end{equation}\\
By $s_\eta(\A)$ we mean the value of the Schur function
corresponding with $\eta$ under substitution of the spectra of $\A$.
An important   fact is that  $dim\ M_\eta (0) = \chi_\eta(e)$, i.e.,
the dimension of the $S_n$ irreducible module corresponding with
$\eta$. See [Kos]or lemma 5.

 It is worth attention  that  formula 9 for  $\hat{d}_\eta(A)$ is the expected contribution of a random subspace $M_\eta (0)$

 to the trace  of a random operator when acts on a whole space $M_\eta$,
 when using an orthogonal form to compute the trace.i e this is the most likely number  one wold evaluate $d_\eta (A )$to A of given spectra
   Hence through  the  monotonicity theorem we have proof  the Merris Watkin monotonicity  conjecture is definitely  supported since the averaging
    values  we compute are
    seemed to be generic   .\\  one can naturally  ask   about  averaging over the  Orthogonal group in place of the Unitarian group.
The answer is quite complicated since $A^{\otimes n} $ usually even
doesn't respect the irreducible modules of the Orthogonal group So
first of all the problem has to be defined more delicately . Next
one has to challenge the difficulties which the Brower algebra
produces . Hopefully it will be treated elsewhere

Now the identity 9 enables us to prove(lemma 7) that for any $S'_n$s
submodule $V \subseteq \otimes^n \mbc^n, \hat{d}_{\chi_V} (A) =
trace (\int {A^u}^{\otimes n} ) \left|_V \right .$

   where $\chi_V$ is the character of $S_n$ action on V.  The trace is defined since such integrals respect $S_n$'s submodule.After this is done ,
    to prove the averaging version of Merris and Watkins  conjecture is a matter of some explicit traces
    computations.(theorem 7)

Sections 2,3,4 can be thought of as preliminaries in representation
theory relevant to our treatment afterwards .  Sections 5,6 deal
with special central element operators on $\otimes^n \mbc^n$ given
by integrations over the Unitarian group.
In Section 7 we prove the monotonicity result.\\
Our treatment doesn't use any explicit integration over the
unitarian group .On the contrary , in   section 8 we bring a formula
for the average of multiplicity-free products of matrix-elements, in
terms of Schur-functions and characters of the symmetric group. The
average is being taken over the unitarian group with respect to its
action by conjugation on the matrix (theorem 9).

\section{The tensor space $\on \quad \quad \mbox{as} \    S_N(U_n) \quad
\quad \mbox{module}$}

Let $v_1, v_2 \cdots v_m$ be an orthonormal basis.  Denote the form by $\langle \ \rangle$.

Given $N$ we extend the form to $\on$ by the formula
\begin{equation}
\langle v_{i_1} \otimes \cdots \otimes v_{i_N} \cdot v_{j_1} \otimes \cdots \otimes v_{j_N} \rangle = \prod_k \langle v_{i_k} v_{j_k} \rangle
\end{equation}\\
Define the Schur action of the symmetric group $S_N$ on $\on$ by the formula
\begin{equation}
\sigma \noo = v_{i_{ \sigma^{-1}(1)}} \otimes \cdots \otimes v_{i_{
\sigma^{-1}(N)}} \quad \quad \sigma \in S_N
\end{equation}
For a given $n \times n$ matrix $A$, define the diagonal action of
$A$ on $\on$ by the formula
\begin{equation}
A \noo = A v_{i_1} \otimes \cdots \otimes A v_{i_N}
\end{equation}
Some times we call the operator coresponding with A" the N Kronecker
power of A" and
denote it by $A^{\bigotimes N}$\\

A monotonic   vector of integers  $\bar{\eta} = \eta_1 \geq \eta_2 \geq \cdots \eta_N$ is called a partition of
$N \ ;  \eta \vdash N$ if $\Sigma \eta_i = N$.\\
For $\eta \vdash N$ define the $\eta$ weight space in $\on \ ; \ \one$ by the formula
\begin{equation}
\on(\eta) = span [ v_{i_1} \otimes \cdots \otimes v_{i_N} ;\#[e : i_e = j ] = \eta_j]
\end{equation}
For a subspace $M$ of $\on$ we define
\begin{equation}
M(\eta) = M \cap \one \ .
\end{equation}
The weight space of $\bar\eta = \eta_1= \eta_2= \eta_3 \cdots$ is
called the zero weight space ; $M(0)$. In this note usually  N=n so
$\eta_i=1$ in the last definition  to zero weight spaces .

Now as it can easily be seen, the action of the symmetric group (11)
and that of the Unitarian group(12)
 on the tensor product spaces  commute with each other
so they respect the isotypical (see section 3) component of one
another . In fact they have an  isotypical component in common  in
their action on $\on \quad \quad $ .So it is enough to treat just
the isotypical component of the symmetric group.\\ because of the
irreducible representations of$ S_N$ indexed by partitions of N it
is reasonable  to denote  the isotypical components of$ S_N(U_n)$
in$\on$ with partitions of N and identify such a  component  with
$V_\eta\otimes M_\eta$ where $V_\eta $ is $S_N$ irreducible and
$M_\eta$is $U_n$ irreducible.\\ The Schur Wyel duality theorem shows
this relationships between $S_n$ and $U_n$ representations on the
tensor product spaces.
\begin{theorem}
As $S_N (U_n) module \on$  isomorphic to $\bigoplus V_\eta\otimes
M_\eta\\The \ sum \ goes  \ over \ \eta \vdash N $ of no more than n
parts

\end {theorem}
Proof can be found in [G W] \\
In the following  sections we treat the action of  the symmetric
group on$\bigotimes^N\mbc^n$ to get more  explicit expression of the
isotypical component in the light of theorem 1 . We shall first give
 some basic information on representation theory .

\section{ representation theory}

Let $G$ be a finite group and $\chi$ an irreducible character of it.  Let $\mbc G$ be the group-algebra of $G$ over $\mbc$.  Define a central element in $\mbc G$ corresponding with $\chi$ by the formula
\begin{equation}
C_\chi = \frac{\chi(e)}{|G|} \sum_{g \epsilon G} \chi(g) \ g
\end{equation}
$C \chi$ is  central  in $\mbc G$ because the characters of $G$ are
conjugacy invariant .More than this.

\begin{lemma}
Given irreducible character $\chi, C_\chi$ is idempotent and \\
$\mbc G \cdot C_\chi$ is the isotypical component of $\on$
corresponding with $\chi$.
\end{lemma}
 To prove this  one  use  the orthogonality relations of the irreducible  characters of a finite group. See [J.K].

For our treatment we formulate it more generally in the next lemma

\begin{lemma}
Let $M$ be a $G's$ module over $\mbc$ and $\chi$ be an irreducible character of $G$.
Then $MC_\chi$ is the isotypical component of $M$ corresponding with $\chi$. \\ .
\end{lemma}

Along this  note  we exchange freely $GL_n (\mbc)$ with the
unitarian group $U_n$ thanks to the next statement.
\begin{lemma}
Each irreducible $GL_n (\mbc)$ module remains irreducible under
reduction to the unitarian group $U_n$ and vice versa, i.e. each
$U_n$'s irreducible module occurs as a reduction from $GL_n (\mbc)$
irreducible module See [GW]  page  94.

\end{lemma}
\textbf{By " isotypical component" }(of a module) we mean  a maximal
submodule   with no  non isomorphic submodules  in it. The
importance of the isotypical component a module  is that
intertwining operators i e operators which commute
 with the action of the group on the module respect it\\
 \textbf{ "complete  reducibility"}is cleared by the next theorem
 \begin {theorem}
 Let G be a compact group (may be finite ) and let V be a G submodule over $\mbc$ than if U is V submodule there exists U'
 a V submodule so that V=$U\bigoplus U'$
\end {theorem}
For proof see [G W]

\section{Partitions, Young diagrams and Young tableaux and the construction of the $U_n \times S_N$ isotypical modules}

For each partition of $\eta$ of $N$, $\eta \vdash N$ one corresponds
irreducible character of $S_N$ $ \chi_\eta$. The explicit
construction of $\chi_\eta$ is not needed for our treatment .
 A good reference  for characters $\chi_\eta$ and their  irreducible module $V_\eta$, is [J.K] or [G W].

Let $\eta \vdash N, \eta = \eta_1 \geq \eta_n \cdots \geq \eta_n$.
For our discussion  any $\eta \vdash N$ is no more than $n$ part.
See Schur Weyl duality( theorem 1)

Now the Young diagram of shape $\eta\vdash N$ is an array of rows of
cells, one under the other, $\eta_i$ cells are in the $i$ row.

 The rows  begin together from the very left to the right, so one
gets an array of columns from left to right as well . See the figure
below.

Given $N$, a Young's tableaux of shape $\eta\vdash N$ is a filling
of the cells of the Young's diagram by the numbers $1, 2 \cdots n$
 so they increase down the column and non-decrease in the rows to the
 right.\\
Next corresponding to partition $\eta \vdash N$ of no more than n parts   we construct a basis for an irreducible  $U_n$ module in $\on$\\
 Given Young tableaux of shape $\eta$, we construct a vector in $\on$ by the next process.
  First, we fix an order on the cells in the Young diagram of shape  $\eta$.  Next, we attach the digits in the
  Young tableaux to the basis vectors of
  $\mbc^n$
  one after another, with respect to the order we had fixed , along a tensor of length $N$ \\Let us take an example .\\
\\
${\young(11,22,3)}\longrightarrow v_1 \otimes v_1 \otimes v_2\otimes v_2\otimes v_3$ .\\
 The order  we  fixed on the cells of the  Young diagram  of shape $\eta = 2 \geq 2\geq 1$
  is that we begin from the upper left down along the rows one after another  ending at the lower right.

Now we act on the tensor we obtained with the central  idempotent  $C \chi_\eta$

Recall $M_\eta$ denotes the $U_n$'s irreducible module corresponding with $\eta \vdash N$.

\begin{theorem}

\begin{itemize}
\item[(i)]Fix an order on the cells of the Young diagram of shape $\eta$.\\
 A basis for $M_\eta$'s copy, in $\on$ is accepted when using the process above over all the Young tableaux of shape $\eta$.
    \item[(ii)] By continuing  the process above in all the orders on the Young diagram of shape $\eta$
    cells,and over all the Young tableau of this shape
    one gets a generating set for the isotypical $U_n$ component(and$ S_N $ as well )of type $\eta$.
\end{itemize}
\end{theorem}
See [K J] or [G W] for proof

\noindent {\bf Corollary} \\
 $dim(M_\eta)$ is the number of Young
tableaux of shape $\eta$.

Let us  define for $\eta \vdash N$ a standard Young tableau as a
filling of the Young diagram of shape $\eta$ with the letters $1 \ 2
\cdots N$ so that they increase down the column and  in the rows to
the right.

\begin{theorem}
$dim (V_\eta)$( $\chi_\eta(e)$), is the number of standard Young
tableaux of shape $\eta$.
\end{theorem}

{\bf Proof.} See [J.K] [G W]

\section{Central $S_N \times U_n$ operators}
\def\eua{E^N_A}
Let $A$ be a $n \times n$ matrix.  Define an operator $\eua$ on $\on$ by the formula
\begin{equation}
\eua = \int (A^u)^{\otimes N} du \quad \mbox

\end{equation}
integrated on the unitarian group with respect to the Haar measure\\
$For\ u \in U_n \ A^u = u A u^* $. \setcounter{lemma}{3}
\begin{lemma}
The operator $\eua$ is $S_N(U_n)$ equivariant.
\end{lemma}
\noindent {\bf Proof.} It is $S_N$ equivariant since for each $u$
$A^{u\bigotimes N}$ is $S_N$ equivariant.  It is $U_n$ equivariant
since the Haar measure is $U_n$ invariant.

Now because $\eua$ is $S_N (U_n)$ commute, it acts scalarly on the
isotypical $S_N \times U_n$ components in $\otimes^N \mbc^n$.

Let $s_\eta(A)$ denote the value of the Schur function  $
s_\eta(\bar{x})$ when  the vector of the eigenvalues of $A$   is
substituted.

\begin{lemma}
$\eua  \left|_{V_\eta \times M_\eta}\right .  = \frac{S_\eta(A)}{S_\eta (I)}$times the identity.  $I$ is the $n \times n$ unit matrix.
\end{lemma}
 \noindent {\bf Proof.} Because $\eua$ is $S_N \times U_n$ equivariant it acts scalarly on the isotypical component, i.e. $V_{\eta} \times M_{\eta}$.
  By definition $s_{\eta}(A) = trace (A) \left|_{M_{\eta}} \right . $ ; $s_{\eta}(I) = dim(M_{\eta})$.  Hence the scalar is $\frac{s_{\eta}(A)}{s_{\eta}(I)}$.

 \noindent {\bf Remark}  $\eua$ depends only on the spectrum of $A$, since the traces of $A$ on $U_n$'s modules
depends only on the spectra of $A$.

It is because the trace of $A$ on each $GL_n(\mbc)$ module depends
just on $A's$ spectra, see [FG].

\section{The traces of $E_{A}^n$ on some subspaces of $\onn$}

Recall that  $M_{\eta}(0)$ is the intersection of $M_\eta$ with the
zero weight space of $\onn$.

 By the last lemma
$trace E_{A}^n   \left|_{M_{\eta}(0)} \right . = dim (M_{\eta}(0))
\frac{s_\eta(A)}{s_\eta (I)}$.
\begin{lemma}
$dim M_\eta(0) = dim (V_\eta)$
\end{lemma}

\noindent {\bf Proof.}  Recall how a basis to $M_\eta$ was
constructed through Young tableaux.( theorem 3) Now the Young
tableaux corresponding with zero weight tensors are those on which
each digit from $1 2 \cdots n$ appears once.  Such a Young tableaux
is standard.  Hence by Theorem 4 the lemma is proved.\\
 .\\
\noindent {\bf Corollary}  $trace \ \eua  \left|_{M_{\eta}(0)}
\right . =  \frac{s_\eta(A)}{s_\eta (I)} \chi_n(e)$.\\
 \noindent{\bf Remark}If one  takes $m\leqslant n$ and  $M_{\eta}(\gamma)$
for $\gamma \vdash m$ which is multiplicity free weight,  the last
lemma remain true (we shall use it in   section 8)

Now we are going to compute  $trace \ \eua  \left|_{M_{\eta}(0)}
\right . $ explicitly by using the standard basis of $\onn (0)$ and
the central idempotent $C_\eta = \frac{\chi_n(e)}{n!} \sum_{\sigma
\in S_n} \chi_{\eta} (\sigma)\sigma $.(Formula 15)  This will be
done along the proof of the next theorem.\\
 We define $\hat{d}_\eta (A) = \int d_\eta A^u du$  integrated  over the unitarian group with respect to the Haar measure on it.

\begin{theorem}
 $trace  E_{A}^n  \left|_{M_{\eta}(0)} \right . = \hat{d}_\eta (A)$
\end{theorem}

\noindent {\bf Proof} $V_\eta \otimes M_\eta = \on \cdot C_\eta$.  Hence  $V_\eta \otimes M_\eta (0) = \on (0) \cdot C_\eta$ .

Hence
\begin{eqnarray}
 trace  E_{A}^n  \left|_{M_{\eta}(0)} \right .& = &  \int du \sum_{\sigma \in S_n} \langle A^u v_{\sigma(1)} \otimes \cdots \otimes v_{\sigma(n)} \cdot C_\eta v_{\sigma(1)} \otimes \cdots v_{\sigma(n)} \rangle \nonumber \\
& = & \int du \ n! \langle A^u v_1 \otimes \cdots \otimes v_n \cdot C_\eta v_1 \otimes \cdots \otimes v_n \rangle
\end{eqnarray}
Here we use the fact that $\sigma C_\eta \sigma^{-1} = C_\eta$.

Now by injection of the explicit expression of $C_\eta$ we come to the last expression
\begin{eqnarray}
\int du \ n! \frac{\chi(e)}{n!} \sum_{\sigma \in S_n} \chi_\eta(\sigma)\langle A^u v_1 \otimes \cdots \otimes v_n \cdot v_{\sigma(1)}
 \otimes \cdots \otimes v_{\sigma(n)} \rangle \nonumber\\ =  \chi_\eta (e) \int  \sum_{\sigma \in S_n} \chi_\eta(\sigma)\prod_i A^u_{i \sigma(i)} du=\chi_\eta (e) \hat{d}_\eta (A)
\end{eqnarray}
At the last step we used the multiplication formula to the scalar
form (Formula 10). The   proof ends  since $\chi_\eta(e)$ is the
multiplicity of $M_\eta(0)$ in $\on$.

To illustrate the last theorem we prove briefly the next theorem of
Merris and Watkins ([M W]theorem 8).

\begin{theorem}
Let $A$ be a $n \times n$ matrix of rank $k$, then $d_\eta (A) = 0$ to any $\eta \vdash n$ but if $\eta$ has no more than $k$ parts.
\end{theorem}

\noindent {\bf Proof.} \\ Recall $s_\eta$ is the Schur function
corresponding  with the partition $\eta$ Now  $s_\eta (A) = 0$ but
if $\eta$ has no more
 than $k$ parts. \\
  Indeed without loss of generality one can assume that A is a diagonal matrix (see the remark at the end of section 5).Now  consider the action of
  $A^{\otimes n}$ on each vector basis corresponding with Young tableaux as in Theorem 3.If $\eta$ has more parts than the rank of A this action
   vanishes identically\\
Now
\begin{equation}
s_\eta (A) = 0 \Rightarrow \hat d_\eta (A)= trac E^n_A \left|
_{V_\eta \otimes M_\eta} \right . = 0
\end{equation}

Since for non-negative definite $A \ d_\eta (A^u) \geq 0$ to any $u \in U_n$ it is to say that $d_\eta (A) = 0$.
To prove the theorem, one can check that the non-negative matrices of rank $k$ are Zarisky dense in the set of matrices of rank $k$.Thanks to
A Goldberger for the remark\\
The next lemma is a clear corollary of theorem(5)and is pivotal to
the proof of the monotonicity theorem
\begin{lemma}
Let $V\subseteq \onn$ be $S_n$'s modules.  Let $\chi_V$ be its character, i.e.
$\chi_V(\sigma) = trace \ \sigma|_{V}$.  Then
\begin{equation}
trace \eua |_{V} = \sum_\sigma \chi_V (\sigma)\int \prod A^u_{i
\sigma i}
\end{equation}
\end{lemma}

\noindent {\bf Proof.}  Because of the complete reducibility of any $S_n$'s module over $\mbc$ and the linearity of the traces and the matrix functions $d_f$
one can deal merely with irreducible modules.

Assume $V$ is irreducible, let us say of type $V_\eta$.  So $V
\subseteq V_\eta \otimes M_\eta$ on which $E^n_A$ acts scalarly,
hence the trace of it depends just on its dimension
 Hence the lemma is proved using  theorem (5)and lemma (6)

\section{The Merris Watkins monotonicity conjecture}
 As  $S_n$ module $\one$, the submodule of tensors of weight $\eta$ , isomorphic to the $S_N$ module induced from the trivial module of the
  Young subgroup corresponding with
  $\eta$ ,$S_{\eta_1}\times S_{\eta_2}\cdots$,. Let say $Young(\eta)$ \\
 i.e $
  \begin{array}{ll} \ 1_{\eta_1} \otimes 1_{\eta_2} \cdots 1_{\eta_n} \end{array}\bigotimes_{\mbc (Young(\eta))}\mbc S_n$.\\ Merris
and Watkins denote the character of this $S_n$ module by $[\eta]$.

Denote by $\succ$ the next relation on the partitions of $n$.

We say that $\eta \succ n'$ if $\sum\limits^j_{i=1} \eta_i \geq
\sum\limits^j_{i=1} {\eta^,}_i$ to any $j$.  Merris and Watkins
conjecture is that $\eta \succ \eta ' \Rightarrow \frac{d_[\eta]
A}{d_[\eta]I} \geq  \frac{d_[\eta '] A}{d_[\eta ']I}$ for any
non-negative definite matrix $A$.one can check that $d_[\eta] I =
\left(
\begin{array}{c} n \\ \eta \end{array} \right)  =
\frac{n!}{\prod_i\eta_i !}$.

 We prove that their conjecture is true when averaging  with respect to conjugacy  relation over $U_n$, the unitarian group.

\begin{theorem}
Under the assumptions above\\
$\eta \succ \eta ' \Rightarrow \left( \begin{array}{c} n \\ \eta \end{array} \right)^{-1} \hat d_[\eta] (A)  \geq \left( \begin{array}{c} n \\ \eta' \end{array} \right)^{-1} \hat d_[\eta ' ] (A) $ to any non-negative definite matrix $A$.
\end{theorem}
The proof includes a series of reductions.\\

First  by lemma (7) for  $\gamma\vdash n$,of n
\begin{eqnarray} \hat d_[\gamma]A^u  = trac E^n_A \ |_{\on (\gamma)} \nonumber\\= \sum\limits \int \langle v_{i_1} \otimes \cdots v_{i_n} A^{u^{\otimes n}} \cdot
v_{i_1} \otimes \cdots \otimes v_{i_n} \rangle \ du
\end{eqnarray}
summed over the tensors of weight$\gamma$\\
Now by using of the product  form for the scaler product formula  over the tensor space (10) we get for (21) the next elagant formula
\begin{eqnarray}
\left( \begin{array}{c} n \\ \gamma \end{array} \right) \int \prod\limits_i \left( A^u \right)^{\gamma_i}_{ii} du
\end{eqnarray}
using the invariancy of the Haar measure to permutations conjugation  we get the next formula
\begin{eqnarray}
  \left( \begin{array}{c} n \\ \gamma \end{array} \right) \int \sum\limits_{\sigma \in S_n} \frac{1}{n!} \prod\limits_i
\left( A^u_{\sigma(i) \sigma(i)} \right)^{\gamma_i}
\end{eqnarray}

 Since A is non negative definite   the diagonal elements  of $A^u$ are non negative
hence the first reduction is to observe that it is enough to prove that for non-negative  $b_1 b_2 \cdots b_n$
\begin{equation}
\sum_{\sigma \in S_n} \prod_i b^{\eta_i}_{\sigma(i)} \geq \sum_{\sigma \in S_n} \prod_i b^{\eta'_i}_{\sigma(i)}
\end{equation}
Next we observe that
\begin{equation}
\sum_{\sigma \in S_n} \prod_i b^{\gamma_i}_{\sigma(i)} = Perm \left(
b^{\gamma_i}_{j} \right) \quad \quad \mbox{(Perm for permanent)}
\end{equation}
and hence then we will prove
\begin{equation}
Perm \left(  b^{\eta_i}_{j} \right)  \geq Perm \left(  b^{\eta'_i}_{j}  \right)
\end{equation}
The next reduction is to assume that $\eta_i = \eta'_i$ to $i \geq 2$.
 We use the invariancy of the permanent we exchange rows  as well as the fact that if $\eta \succ \eta'$ there exists a path of partitions between
 $\eta$ to $\eta'$ monotonic with respect to $\succ$ so that each consecutive partition along the path differ just on two  parts.

Next we develop the  last permanent in respect to their   first two
rows.We show that in each summand in the computation we get the
desired inequality
   so  the last reduction is to prove the next lemma.

\begin{lemma}
Let $A,B$ be positive and for a given $n$ denote \\$\varphi_{A,B}(x) = per \left( \begin{array}{ll} A^x & B^x \\
 A^{n-x} & B^{n-x} \end{array} \right)$
then for $x \geq \frac{n}{2} \; \; \varphi_{AB}(x)$ increases with $x$\\
\noindent {\bf Proof.} \\
 We consider $\frac{d}{dx} \varphi_{AB}(x)$
$$
\frac{d}{dx} \varphi_{AB}(x)  =  \frac{d}{dx} (A^{x} B^{n-x} + A^{n-x} B^x) = A^x B^{n-x} \log \frac{A}{B} + A^{n-x} B^x \log \frac{B}{A}
$$
$$
= [A^x B^{n-x} - B^x A^{n-x}] \log \frac{A}{B}  = A^x B^{n-x} \left( 1 - \left( \frac{B}{A} \right)^x \left( \frac{A}{B} \right)^{n-x}  \right) \log \frac{A}{B}
$$
now if $A \geq B$ then $\log \left( \frac{A}{B} \right) > 0$ and $\left( \frac{B}{A} \right)^x \left( \frac{A}{B} \right)^{n-x} \leq 1$
 since $x \geq \frac{n}{2}$.\\
  Hence the derivation is positive .  The other case is treated the same way.
\end{lemma}
This ends the proof of  theorem (7).

As an example of using the monotonicity theorem,  we give a proof to the theorem of James and Lieback [JL] on the dominancy of the permanents among the immanents of no more than two parts, but in average.

\begin{theorem}
Let $\eta_1 \geq \eta_2$  be a partition  of two parts then \\$\frac{\hat{d}_\eta(A)}{d_\eta(I)} \leq per (A)$ for any non-negative definite $n \times n$ matrix.
\end{theorem}
\noindent {\bf Proof.} As it is well known [JL], [JK],   $\chi_\eta
= [\eta] - [\eta']$ when $\eta' = \eta'_1 \geq \eta'_2$ so that
$\eta'_1 = \eta_1 + 1$ so
\begin{equation}
\frac{\hat{d}_[\eta](A)}{d_[\eta](I)} = \frac{\hat{d}_[\eta'](A)+
\hat{d}_\eta(A)}{\hat{d}_[\eta'](I)+d_\eta(I)} =
\frac{\frac{\hat {d}_[{\eta'}](A)}{d_[\eta'](I)} d_[\eta'](I) +
\frac{\hat{d}_\eta(A)}{d_\eta(I)} d_\eta(I)  }     {d_[\hat{\eta'}](I)
+ d_\eta (I)  }
\end{equation}
The last expression is a convex sum of
$\frac{\hat{d}_{[\eta']}(A)}{d_[\eta'](I)}$ and
$\frac{\hat{d}_\eta(A)}{d_\eta(I)}$.  Hence it is greater than the
minimum. Assume $\frac{\hat{d}_\eta(A)}{d_\eta (I)} >
per(A)\geq\frac{\hat{d}_{[\eta']}A}{d_[\eta'](I)}$ one gets
$\frac{\hat{d}[\eta](A)}{d_[\eta](I)} > \frac{\hat{d}_[\eta'](A)}{d_[\eta'](I)}$ and this contradicts the  monotonicity we have  just proved.\\
Remark on generalizations\\
One can check  that under the next definition of \\ $\hat d_\eta(A)\ for\ \eta\vdash\ m \ \leqslant n$  theorem 7 remain  true\\
The definition is given in the next formula\\

For  $m\leqq n$  and $ \eta \vdash m $
 \begin{equation}
 \hat{d}_\eta(A)=\int\sum_{\sigma\in S_m}\chi_\eta(\sigma) \prod_{i\leqq m}A^u_{i\sigma(i)} du
\end{equation}
integrated over the unitarian group $U_n$\\Especialy lemma 6 remain
true for $m\leqq n$.\\
With this last remark we are coming to the last section  dealing
with some explicit expressions for the integration of matrix
monomials which have occurred  over the note
\section{$ U_n$ Invariant  matrix's elements' products}
Let A be  $n\times n$ complex matrix and $\lambda_1, \lambda_2...\lambda_n$ its eigenvalues
For given $m\leq n$ let $\sigma\in S_m$ .
For $\gamma\vdash m$ let $s_\gamma$ be its corresponding Schur function\\
Let us define $I(A,\sigma)=\int\prod A^u_{i\sigma i}du$ integrated over the unitarian group with respect to the Haar measure\\
\\
\begin{theorem}
$I(A,\sigma)=\sum_{\eta\vdash m}\frac
{s_\eta(\lambda_1...\lambda_n)}{s_\eta(1,1...1)}\frac{\chi_\eta(e)}{m!}\chi_\eta(\sigma)$
\end{theorem}
 We give several example before the proof\\

1)Let A=I the unit $n\times n$ matrix\\
Using the formula one get\\
$I(I,\sigma)=\sum_{\eta\vdash m}\frac{\chi_{\eta(e)}\chi_\eta(\sigma)}{m!}=\delta_{e,\sigma}$\\
Indeed this is the second  orthogonal relatione of the characters of the symmetric group\\
\\
2)Let A be $2\times 2$ matrix and $\lambda_1,\lambda_2$ its eigenvalues \\
The  two partitions of 2 are the trivial (2) and the only non trivial(1,1)\\
 Now \\$s_{(2)}=\lambda_1^2+\lambda_2^2 +\lambda_1\lambda_2$\\
 $s_{(1,1)}= \lambda_1\lambda_2$ \\
 hence $I(A,\sigma=1,1)=\frac{1}{2}\frac{\lambda_1^2+\lambda_2^2+\lambda_1\lambda_2}{3}+\frac{\lambda_1\lambda_2}{2}$\\
 $I(A,\sigma=(12))=\frac{1}{2}\frac{\lambda_1^2+\lambda_2^2+\lambda_1\lambda_2}{3}-\frac{\lambda_1\lambda_2}{2}$\
 Hence the determinant ie  $I(A,\sigma=1,1)-I(A,\sigma=(12))=\lambda_1\lambda_2$ \\

3)Let A be a matrix of rank 1 than
 $I(A,\sigma)=\frac {s_m(\lambda_1,0,0...0)}{s_m(1,1,1..)}$ where $\lambda_1$ is the only non zero eigenvalue of A. Since only the one part
 partition Schur function supports a matrix of rank 1(Recall  theorem 6)
  \\
 Hence the integral doesn't depend on $\sigma$\\
 We turn to prove the theorem\\
 by lemma 4and Schur Weyl duality since
 $E_A^m$ is $GLn(\mbc)$,equivariant   one can write
$E_A^m=\sum a_\sigma \sigma$ summed over the symmetric group \\
Let us compute the coefficients $a_\sigma$ \\
By lemma 5  $\eua  \left|_{V_\eta \times M_\eta}\right .  =
\frac{S_\eta(A)}{S_\eta (I)}$ times the identity  . On the other
hand by lemma 1 $C_\eta$ acts as the unit on $V_\eta \times M_\eta$
and is vanished on the other isotypical components so using formula
15 to $C_\gamma $ we get
\begin{equation}
 E^m_A= \sum\frac{s_\eta(A)}{s_\eta (I)}\sum \frac{\chi_\eta(e)}{m!}\sum \chi_\eta(\sigma)\sigma
 \end{equation}
summed over $S_m$ and over all partitions of m  . Now for $ m \leqq
n$ we define the zero weight space of $\otimes^m\mbc^n$ to be  the
span of
 all the tensors of type
 $(v_{\sigma (1)}\otimes...\otimes v_{\sigma (N)})$
 for $\sigma\epsilon S_m$\\
  For $\eta \vdash m \ M_\eta(0)$ be the intersection of $M_\eta$
 with the zero weight space\\
   For $\gamma\vdash m$ we compute $trace  E_{A}^m|  _{M_\gamma(0)\bigotimes V_\gamma}$ in two
 ways\\
 The first one is to substitute $\gamma$ in formula 29 (it means $\sigma\Rightarrow\chi_\gamma(\sigma)$) and multiply by $dim(M_\gamma(0))$. \\
 Now $dim(M_\gamma(0))=\chi_\gamma(e)$ by the generalization remark
 after theorem 7.
 On the other hand we compute $trace E_{A}^m|  _{M_\gamma(0)\bigotimes V_\gamma}$ explicitly using the basis of the zero weight space

 \begin{equation}
\sum \int< A^{u \otimes
m}(v_{\sigma(1)}\otimes....v_{\sigma(m)})C_\gamma(v_{\sigma(1)}\otimes....v_{\sigma(m)})>du
\end{equation}
As in theorem 5 one can reduce formula 30 to
 \begin{equation}
  \chi_\gamma(e)\int\sum
<A^{u\otimes m}(v_1\otimes...v_m)(v_{\sigma (1)}...\otimes v_{\sigma
(m)})>\chi_\gamma(\sigma)du
\end {equation}
summed over $S_m$
Now by use of the product rule of the scaler product (formula 10) in $\otimes^m\mbc^n$ one get\\
$\chi_\gamma(e)\sum\int \prod_i A^u_{i\sigma(i)}\chi_\gamma(\sigma
)du$ summed over $S_m$

 Now one get theorem 9 by equating of coefficients of $\chi_\gamma
$ in the two ways of the trace computations \noindent \\{\bf Remark}
Kavin Coulembier  ,in a lecture  given in Decin's conference at
August 2011, pointed out  the similarity   of such integrals over
the unitarian group to those integrals in [Co Sn]  at least by using
the Schur Weyl duality  theorem , indeed they  also used Schur Wyel
duality in the case of  Orthogonal (Symplectic ) groups but they
reduced their attention to groups' matrices .See our remark in the
introduction about averaging over the Orthogonal groups problems

\section*{References}
\begin{itemize}

\item[{[Kos]}] B. Kostart.  Immanent's inequalities and o weight spaces.  Journal of the AMS 8, (1995), 181-186.
\item[{[MW]}] R. Merris and W. Watkins.  Inequalities and identities for generalized matrix functions.  Linear Algebra and its Applications, 64:(1985), 223-242.
\item[{[Ch.Ha Mit]}]M. Christandl; A. Harrow; G. Mitchison Non-zero Kronecker coefficients and consequences for spectra.  Communication in Math Physics, 270(3) (2007), 575-585.
\item[{[GF]}] A. Goldberger, A. Frumkin.  On the distribution of the spectrum of the sum of two Hermitian or real symmetric matrices.  Advances in Appl. Math. 37 (2006), 268-286.
\item[{[WG]}] R. Goodman, N.R. Wallach.  Representations and invariants of the classical groups. (2003).
\item[{[P 1]}] Thomas A. Pate.  Descending chains of immanent.  Linear Algebra and its Applications. 162-164 (1992), 639-650.
\item[{[P 2]}] Thomas A. Pate  Row appending maps.$\psi$functions and immanent inequalities for Hermitian positive semi definit matrices
                               Proceeding of the London math society 1998 76 307-358
 \item[{[P 3]}]Thomas A. Pate. Psi functions,irreducible characters and Merris and Wotkins conjecture lin multi lin algebra 35 1993 195-213
\item[{[Lieb]}] E.H.Lieb  Proof of some conjecture on permanents  Journal of math and mechanics 16(1966)127-134
\item[{[JL]}] G.D. James and M. Liebeck.  Permanents and immanents of Hermitian matrices.  proceeding of the London Math. Soc. 55(3), (1987) 223-242.
\item[{[Li]}]D Littlewood . The theory of group characters, (1950).
\item[{[Co Sn]}]B. Collins .P .Sniady .Integration with respect to the Haar measure on unitary, orthogonal and symplectic group  Commun. Math. Phys. 264, 773–795 (2006)
\end{itemize}

\end{document}